%% file: stb.tex
\documentclass[11pt,a4paper,final,leqno,notitlepage]{amsart}
\usepackage[dvips,final]{epsfig}
\usepackage{epsf, graphics, latexsym, eucal, pstcol,
amsfonts, amstext, amssymb, amsmath, amsthm, amsbsy, amsxtra, amscd }

\theoremstyle{plain}

\newtheorem{lstheorem}{Theorem}[section] %% theorems, labelled 1.1, 1.2,...
\newtheorem*{corollary}{Corollary}    
\newtheorem{lsproposition}{Proposition}[section]

\newtheorem*{proposition}{Proposition}

\theoremstyle{definition}

\newtheorem{lsdefinition}{Definition}[section]

\theoremstyle{remark}

\newtheorem*{remark}{Remark}

\begin{document}

\hoffset = 0pt
\voffset = 0pt
\topmargin = 0pt
\headheight = 0pt
\textwidth = 355pt

\setlength{\leftmargin}{0.1cm}
\setlength{\rightmargin}{0.1cm}

\renewcommand{\abstractname}{Abstract}

\title{On the stability of dual scattering channel schemes}

\author{Steffen Hein}
\input{abs.tex}   %abstract

\maketitle
\pagestyle{myheadings}
\markboth{{\normalsize \textsc{Steffen Hein}}}
{{\normalsize \textsc{Stability of DSC Schemes}}}
\normalsize
\vspace{-0.9cm}
\input{int.tex}   %introduction 
\input{sec2.tex}  %section 2
\input{sec3.tex}  %section 3
\input{ccl.tex}   %conclusions

\input{ref.tex}

\vspace{5pt} 
\hspace{-16pt}
{\Small
\address{\textsc{Spinner RF Lab, Aiblinger Str.30,
DE-83620 Westerham, Germany}}
\newline\!E-mail address:\; \email{s.hein@spinner.de}
}
\end{document}

%% file: abs.tex
\begin{abstract}
Dual scattering channel (DSC) schemes generalize \textsc{Johns'}
%\linebreak
TLM algorithm in replacing transmission lines with abstract scattering
\linebreak
channels in terms of paired distributions. A well known merit of TLM
\linebreak
schemes is unconditional stability, a property that is commonly drawn
%\linebreak 
upon the passivity of linear transmission line networks. So the question
%\linebreak
arises, if DSC algorithms remain stable in a neat sense. It is shown
%\linebreak
that a large class of \emph{$\alpha$-passive} processes are in fact
unconditionally stable.
%\linebreak
The ana\-ly\-sis applies to TLM and DSC schemes alike and includes
\linebreak 
non-linear situations.

\vspace{.2cm}
\hspace{-.55cm}
%Westerham in summer 2004
Westerham on June 8, 2004
%Westerham on \today
\end{abstract}

%% file: int.tex
\section{Introduction}\label{S:sec1}
Dual scattering channel (DSC) schemes result from an incisive revision
of the transmission line matrix (TLM) numerical method.
The latter has originally been introduced by \textsc{P.B.Johns} and
coworkers in the early 1970s ~\cite{JoBe} and has since been subject
to assiduous study and publication ${\text{[Tlm1-3]}}$.
The TLM method is today commonplace in scientific computing and largely 
applied to the numerical solution of \textsc{Maxwell}'s equations ~\cite{Hoe}
but also to manifold wave propagation, transport, and diffusion phenomena;
we can here refer to the monographs of \textsc{Christopoulos}
~\cite{Ch} and \textsc{de Cogan} ~\cite{dC}. Also, \textsc{Rebel} ~\cite{Re}
gives a fairly complete survey over the state of the art of TLM by the
year 2000.

DSC schemes are generalized TLM methods in arising from a twofold
abstraction ~\cite{He1}.
Firstly, the scattering channel concept underlying TLM is redefined in
terms of paired distributions. 
(Characteristic impedances are thus neither needed, nor in general
defined, e.g.)
In the second place, non trivial cell interface scattering is admitted
during the connection step of iteration, thus taking advantage of the
intrinsic duality in the connection-reflection cycle of the algorithm.

The extended framework bypasses a set of modeling limitations induced by
transmission lines (\,and discussed in more detail in ~\cite[section2]{He1}\,)
while it preserves the main advantages of the TLM method.
In particular, the convolution type updating scheme and Johns' two step
connection-reflection cycle are essentially retained.
The question arises, however, if DSC schemes remain \emph{unconditionally
stable} in a well-defined sense.

TLM schemes are unconditionally stable in that they are equivalent
to passive linear transmission line network models ~\cite{Jo2}.
A concept that applies to DSC and TLM algorithms alike
is $\,\alpha$-\emph{passivity} which paraphrases contraction properties
with respect to a non-negative (de)limiting functional; \linebreak
cf. section ~\ref{S:sec2}.
In the TLM context, ${\,\alpha\,}$ is essentially the sum over the
squared incident and outgoing transmission line voltages. The sum is
contracted at each scattering event due to energy conservation
(or loss).

In a like manner $\,\alpha$-passivity characterizes the reflection
and connection maps of an unconditionally stable DSC algorithm,
but the 'energy' functional $\,\alpha$ (which may in fact measure
\emph{any} conserved quantity) needs not to be a quadratic form.
For instance, energy may be linearly related to temperature and
quadratically to particle velocity, within the same algorithm.
The following approach is sufficiently general to apply to such
and other, even non-linear situations.

%% file: sec2.tex
\section{Stability}\label{S:sec2}
Algorithm stability prevents the computational process from piling up
to infinity (it does not yet imply convergence, or ~\emph{consistence}
of the algorithm, of course).
TLM models are \emph{unconditionally\/} stable in that they are
equivalent to passive linear transmission line networks ~\cite{Jo2}.
DSC schemes, in not using lines, need a more general characterization
which here is given in terms of $\alpha$-passive causal functions.

Let ${\mathcal{L}\/}$ be a real or complex linear space and $I$ 
a totally ordered set 
(\,e.g. $I \in \{\mathbb{N}\,, \mathbb{Z}\,, \mathbb{R}\}\,$;
intervals are then naturally defined in $\,I\,$ by the order relation.
We commonly think with $\,I\,$ of a discrete or continuous time domain.)
Also, let ${\mathcal{E} \subset \mathcal{L}^{I}\/}$ be a set
of functions such that ${f\in\mathcal{E}\,}$ implies
$\chi_{_{s\leqslant t}}(s) f(s)\in\mathcal{E}\,$,
for every ${t\in I\,}$\,, where $\chi_{_P}(s)$ denotes the
characteristic function of property $\,P\,$\, 
(which is 1 if $\,s\,$ shares that property and 0 else).

\begin{lsdefinition}\label{D.2.1}
A function
$F:\mathcal{E}\to\mathcal{E}\,$ is called \emph{causal\/},
iff for every ${f\in\mathcal{E}\,}$\\ and ${t \in I\,}$
\vspace{-.25cm}
\begin{equation}\centering\notag %\label{2.0}
F\,f\,(\,t\,)\;=\;F\;[\;\chi_{_{s\leqslant t}}f\,(\,s\,)\;]\,(\,t\,)\quad .
\end{equation}
Such functions are also called (causal) \emph{propagators}.

\begin{remark}
In some respect, causal functions generalize lower triangular matrices
or integral operators such as
\vspace{-.25cm}
\begin{equation}\centering\notag %\label{2.1}
F\,f\,(\,t\,) \; = \; \int_{-\infty}^{t}K\,(\,t\,-\,s\,)\,f\,(\,s\,)\; ds
\quad
\end{equation}
with a Green's function kernel $K$, e.g.
Note that in general $\,F:f\mapsto F f\,$ needs not to be linear.

Typically $\,f\,$ represents a state evolving in time
(\,i.e. a \emph{process}\,).
Then causality of $\,F\,$ means that $\,F\,f\,(\,t\,)\,$
depends on the history of $\,f\,$ only up to \emph{present}\,
time $\,t\,$.
\end{remark}
\end{lsdefinition}
\hspace{-.53cm}
The proof of the following is easy and left to the reader.
\begin{proposition}\label{P.2.1} \ \\
For every ${\,t\in I\,}$;
${\,f,\,g\,\in\mathcal{E}\,}$ and causal functions
${F,\,G\,:\,\mathcal{E}\to\mathcal{E}}$
\begin{itemize}
\item[(i)]
$F f\,(\,s\,)\;=\;F\;[\;\chi_{_{r\leqslant t}}(\,r\,)\,f\,(\,r\,)
\;]\,(\,s\,)\,$,
for every $\,s\leqslant\,t$ \;.
\vspace{4pt}
\item[(ii)]
If $\,f(s)\;=\;g(s)\,$ for \emph{every} $\,s\leqslant\,t\;$, \\
\vspace{-0.45cm}
\begin{equation}\centering\notag %\label{2.2}
\text{then} \quad F f\,(s)\;=\;F g\,(s)\;\,\text{for every\,}
\,s\leqslant\,t\;.
\end{equation}
\end{itemize}
\begin{itemize}
\item[(iii)]
The product of causal functions
\begin{equation}\centering\notag %\label{2.3}
\begin{aligned}
FG \, : \, \mathcal{E}\;&\longrightarrow\; \mathcal{E} \\[-2pt]
f\;&\longmapsto\;F\,G\,f\,\;:\,=\;F\,[\,G\,[\,f\,]\,]
\end{aligned}
\end{equation}
is again causal. In fact, if $\mathcal{E}$ is a linear space,
then the causal functions over $\mathcal{E}$ form an algebra.
\end{itemize}
\end{proposition}

\hspace{-.55cm}
Let ${\|...\|\,}$ be a norm on $\mathcal{L}\,$ and 
$\alpha \in \mathbb{R}^{\,\mathcal{L}}\,$ a continuous non-negative real
functional on $\,\mathcal{L}$.
\begin{lsdefinition}\label{D.2.2}
\begin{itemize}
\item[]
\item[(i)]
We call a process ${g\,:\,I\to\mathcal{L}}$ \emph{stable},
iff $\,g\,$ is uniformly bounded on $I\,$,
i.e.\! iff there exists ${\,b\in\mathbb{R}_{+}\,}$ such that 
${\|\,g\,(t)\,\|\,<\,b\,}$ for every ${t\in I\,}$.
\item[(ii)]
The functional ${\alpha\,:\,\mathcal{L}\to\,\mathbb{R}\,}$
is named a (\emph{de})\emph{limiting functional\,},
iff there exist any non-negative real constants 
$\,a,b,c\,$ such that
\begin{equation}\centering\label{2.4}
\|\,z\,\|\;\leqslant\;a\,+\,b\,(\alpha\,(\,z\,)\,)^{c}\quad,
\end{equation}
for every ${\,z\in\mathcal{L}\,}$.
Then obviously ${\,b, c > 0\,}$, and we say also that $\alpha\,$ 
is \emph{minimal increasing\,} (\emph{in any order\,})
\emph{not lesser than} ${\,1/c\,}$.
\end{itemize}
\end{lsdefinition}

Let $\,\mu\,$ be a measure on $\,I\,$ such that intervals are
$\,\mu$-measurable sets.
Functions on $\,I\,$ are henceforth read \emph{modulo $\mu$}
(viz. as equivalence classes of functions that differ at most on sets of
$\mu$-measure zero). Also, let
${\,\alpha\in\mathbb{R}^{\,\mathcal{L}}\,}$
be a delimiting functional on $\mathcal{L}\,$
(i.e. one that is increasing not lesser than any positive order ${\,1/c\,}$),
and assume that ${\,\alpha\circ f\,}$
is $\,\mu-$summable over finite intervals in $\,I\,$
for every ${\,f\in\mathcal{E}\,}$.
The latter is for instance the case if 
$\,\alpha\,(\,z\,)=\|\,z\,\|^{\,p}$
for any real ${\,p\geqslant 1\,}$
and ${\,\mathcal{E}\,\subset\,L^{\,p}(\,I,\,\mathcal{L}\,)\,}$,
which is the metric completion of
\vspace{-.2cm}
\begin{equation}\centering\notag %\label{2.5}
\{\, f \, \in \mathcal{L}^{I}\, \mid \,
(\,\int\nolimits_{I}\|\,f\,\|^{p} \;d\mu \;)^{\;\frac{1}{p}}
\; \leqslant \; \infty \} \quad,
\end{equation}
i.e. $\mathcal{E}\;$ is a subset of the Banach space with norm
$\,\|\,f\,\|_{p}\;:\,=\;
(\,\int\nolimits_{I}\|\,f\,\|^{p}\;d\mu\;)^{\,\frac{1}{p}}\,$.

\begin{lsdefinition}\label{D.2.3}
A causal function $\,F:\mathcal{E}\to\mathcal{E}\,$
is called $\alpha$-\emph{passive\/}, iff
\begin{equation}\centering\label{2.6}
\int\nolimits_{s<t}\alpha\,(\,F f\,(\,s\,)\,)\;d\mu(s)\;\leqslant\;
\int\nolimits_{s<t}\alpha\,(\,f\,(\,s\,)\,)\;d\mu(s)\quad,
\end{equation}
for every ${\,f\in\mathcal{E}\,}$ and ${\,t\in I\/}$.
\begin{remark}
If $\,\alpha\,=\,\|...\|^{p}$ for any real $\,p\geqslant 1\,$ and 
${\,\|\,f\,\|^{p}\,}$ is $\,\mu$-summable over $\,I\,$,
i.e. $\,f\in L^{\,p}(\,I,\,\mathcal{L}\,)\,$,
then \eqref{2.6} clearly implies 
$\,F f\, \in L^{\,p}(\,I,\,\mathcal{L}\,)\,$ and
\begin{equation}\centering\notag %\label{2.7}
\|\,F f\,\|_{_{p}} \;\leqslant\;
\|\,f\,\|_{_{p}} \quad.
\end{equation}
Hence, every ${\|...\|^p}$-passive causal function $\,F\,$ 
defines a \emph{contraction} operator on
${\mathcal{E}\,\cap\,L^{\,p}(\,I,\,\mathcal{L}\,)}$.
\end{remark}
\end{lsdefinition}

Assume, furthermore, that ${\tau\in\mathbb{R}_{+}\,}$ and let on
${\,I\;:\,=\;\mathbb{R}\,}$ the measure $\mu$ be concentrated in
${\{\,k\tau\mid k\in\mathbb{Z}\,\}\,}$ with uniform weight
${\mu\,(\,\{k\tau\}\,)\,=\,\tau\,}$, ${k\in\mathbb{Z}\,}$. 
Alternatively, let $I\;:\,=\;\{\,k\tau\mid k\in\mathbb{Z}\,\}\,$
with 'the same' measure $\,\mu\,$.                               \linebreak 
(Virtually we deal of course with that
discrete situation, even in working on the real axis
with functions that are constant over intervals
${[\,k\tau, (k+1)\,\tau)\,}$; ${k\in\mathbb{Z}\,}$.
Indeed we retain the integral formalism for simplicity, and
the reader may optionally re-write the following integrals as sums.)
For every $f\in\mathcal{E}$ let ${f(t-\tau)}$ $\in\mathcal{E}$,
i.e. $\,\mathcal{E}\,$ is closed under time shifts by negative
integer multiples of $\,\tau$.
Then, for arbitrary ${N\in\mathbb{N}_{+}\,}$ and \emph{exciting function\,}
$e\in\mathcal{E}\subset\mathcal{L}^{I}\,$ with support on
${[\,0,\,N\tau\, )\subset I\,}$, the following holds

\begin{lstheorem}{\label{T.2.1}}
\emph{[\,Stability of the iterated passive causal process\,]}
\newline
For any $\alpha$-passive causal function
$\,F:\mathcal{E}\to\mathcal{E}\,$ and
${\,e \in \mathcal{E}\,}$ as stated,
\newline if ${g\in\mathcal{L}^{I}\,}$ is a process such that
${g(t)\,\equiv\,0\,}$ for $\,t\leqslant 0\,$ and
recursively for \\
${\,t\,=\,n\tau;\,n\in\mathbb{N}\,}$
\vspace{-.25cm}
\begin{equation}\centering\label{2.8}
g\,(\, t \, + \tau \,) \;=\; \; F\,[\,e\,+\,g\,]\,(\,t\,) \quad,
\end{equation}
then $\,g\,$ is \emph{uniquely defined}
\textnormal{(}\emph{modulo $\mu$}\textnormal{)}
and for every $t\geqslant N\tau\,$ holds
\vspace{-.10cm}
\begin{equation}\centering\label{2.9}
\|\,g\,(\,t\,)\,\|\;\;\leqslant\;\;a\;+\;
(\;\;\frac{b}{\tau}\;\int_{[\,0,\,N\tau\,)}\alpha\,(\,e\,+\,g\,) -
\alpha\,(\,g\,)\;\;\;d\mu\;\;)^{\;c}\quad\qquad
\end{equation}
with every constants $a,b,c\,$ that satisfy \eqref{2.4}.
Hence $g\,$ is \emph{stable}.
\end{lstheorem}

\begin{remark}
A process $\,g\in\mathcal{L}^{I}\,$ which is recursively generated
according to \eqref{2.8} by iteration of an $\,\alpha$-passive causal
function $\,F$ is called an $\alpha$-\emph{passive process}.

The theorem ensures thus that for any excitation of finite duration 
(and with no further restrictions) the $\alpha$-passive process
is necessarily stable and in this sense \emph{unconditionally} stable.

Note that \emph{existence} of such a process $g$ is not apriori guaranteed,
since this obviously depends on the condition that with
${s_{0}\,:\,=\,e\in\mathcal{E}\,}$ recursively also the functions
\begin{equation}\centering\label{2.10}
s_{n}(t)\,:\,=\,e(t)+F[\,s_{n-1}\,]\,(\,t-\tau\,)\;\in\;\mathcal{E}\;,
\end{equation}
for $\,0\,<\,n\,<\,N\,$ 
- ~which eventually has to be checked.
\newline Clearly, conditions \eqref{2.10} are always true (and hence $\,g\,$
exists) if $\,\mathcal{E}\,$ is a linear space.
We do not universally premise this, in order to potentially apply the theorem
to non-linear situations, where conditions \eqref{2.10} may only be
satisfied for sufficiently small excitations $e\,$.
\end{remark}

\begin{corollary}
\begin{itemize}
\item[]
\item[(i)]
In the special case $\,\alpha\,=\,\|...\|\,$ estimates \eqref{2.4} holds
with ${a=0}$, ${b=1}$, ${c=1}$.
Then the triangle inequality applies to the
integrand of \eqref{2.9} and validates the bound
\begin{equation}\centering\notag %\label{2.11}
\,\|\,g\,(\,t\,)\,\| \;\; \leqslant \;\;
\frac{1}{\tau}\;\int_{[\,0,\,N\tau\,)}\|\,e\,\|
\;\;d\mu \quad.
\end{equation}
\item[(ii)]
If $\,N=1\,$, i.e. ${\,e(t)\,}$ is a \textsc{Dirac} excitation concentrated
on ${\,[\,0\,,\tau\,)\,}$ \,\textnormal{(}where ${\,g=0\,}$\textnormal{)},
then \eqref{2.9} reads simply
\begin{equation}\centering\notag %\label{2.12}
\|\,g\,(\,t\,)\,\|\;\leqslant\;a\,+\,
(\;\frac{b}{\tau}\;\int_{[\,0,\tau\,)}\alpha\,(\,e\,)
\;\;d\mu\;\;)^{\,c} \quad,
\end{equation}
provided that $\,\alpha\,(\,0\,)\,=\,0\,$ \textnormal{(}
which is the normal case \textnormal{)}.
\end{itemize}
\end{corollary}

\begin{proof}
Clearly, ${\,g\in \mathcal{L}^{I}\,}$ is uniquely defined by the
given recurrence relations, since ${\,e\,+\,g\,}$ (\,and hence $\,g\,$\,)
at the right hand side of \eqref{2.8} is evaluated only up to time
$\,t\,=\,n\,\tau\,$, in virtue of the causality of $\,F\,$.
\newline Furthermore, if
$\,N\leqslant n\,$ and $a\,\leqslant\,\|\,g\,(\,n\tau\,)\,\|\,$
with any $a\,$ satisfying \eqref{2.4}, then
\newline with pertinent $b,c\,$ that satisfy \eqref{2.4}
\begin{equation}\centering\notag %\label{2.13}
\begin{aligned}
0\;\leqslant\;(\,\frac{1}{b}\,
&(\,\|\,g\,(\,n \tau\,)\,\|\,-\,a\,)\,)^{\,1/c}\;
\leqslant\;\alpha\,(\,g\,(\,n \tau\,)\,)\\
&=\;\frac{1}{\tau}\;(\;\int_{s < (n+1)\,\tau}\,\alpha\,(\,g\,)\;d\mu(s)
\; - \;\int_{s < n \tau}\,\alpha\,(\,g\,)\;d\mu(s)\;)\\
&=\;\frac{1}{\tau}\;(\;\int_{s < n \tau}\,\alpha\,(\,F[\,e\,+\,g\,]\,)\;d\mu(s)
\; - \; \int_{s < n \tau}\,\alpha\,(\,g\,)\;d\mu(s)\;)\\
&\;\;\rotatebox{90}{$\rightsquigarrow$}
\text{\small{\;recursion formula \eqref{2.8}}}\\
&\leqslant\;\frac{1}{\tau}\;(\;\int_{s<n\tau}\,\alpha\,(\,e\,+\,g\,)\;d\mu(s)
\;-\;\int_{s<n\tau}\,\alpha\,(\,g\,)\;d\mu(s)\;)\\
&\;\; \rotatebox{90}{$\rightsquigarrow$}
\text{\small{\;since F is passive}}\\
&=\;\frac{1}{\tau}\;\int_{[\,0,\,N\tau\,)}\,\alpha\,(\,e\,+\,g\,)
\,-\,\alpha\,(\,g\,)\;\;\;d\mu(s)\;\\
&\;\; \rotatebox{90}{$\rightsquigarrow$}
\text{\small{\;since $\,e(\,t\,)\equiv 0\,$\; if
$\,t \notin [\,0,\,N\tau)\,$}}\quad.
\end{aligned}
\end{equation}
Thus, estimates \eqref{2.9} holds true in the case
$a\,\leqslant\,\|\,g\,(\,n\tau\,)\,\|\,$ and trivially otherwise.
It follows that ${\,\|\,g\,\|\,}$ is uniformly bounded on 
${\,I \smallsetminus [\,0,\,N\tau)\,}$, hence also on $\,I\,$, 
that is to say $\,g\,$ is stable.
\end{proof}

%% file: sec3.tex
\section{DSC Processes}\label{S:sec3}
In this section DSC schemes are represented as paired $\,\alpha$-passive
processes such as are dealt with in Theorem ~\ref{T.2.1} 
and which hence are stable.

Just as the TLM algorithm, DSC schemes operate on a space
$\,\mathcal{P}\,$ of \emph{propagating fields},
which is a product of (real or complex) normed linear spaces
\begin{equation}\centering\label{3.1}
\mathcal{P}\;=\;\;\mathcal{P}_{in}\,\times\,\mathcal{P}_{out}\quad,
\end{equation} 
cf. ~\cite{He1}.
The two factors are named (somewhat off-hand) the \emph{incident} and
\emph{outgoing} subspace of ${\,\mathcal{P}}$. They are isomorphic in
that there is canonical involutary isomorphism of normed linear space
\begin{equation}\centering\label{3.2}
\begin{aligned}
nb\;:\;\mathcal{P}\;&\to\;\mathcal{P}\\
z\;=(\,z_{in}\,,\,z_{out}\,)\;&\mapsto\;
(\, z_{out}\,,\,z_{in}\,)\;=\,:\; nb\,(\,z\,)\;,
\end{aligned}
\end{equation}
which is commonly called the \emph{node-boundary} map.
Hence, there exists a space
${(\,\mathcal{L}\,,\,\|...\|\,)\,}$
such that 
\begin{equation}\centering\label{3.3}
\mathcal{P}_{in}\;\cong\;\mathcal{P}_{out} 
\;\cong\;(\, \mathcal{L}\,, \,\|...\|\,)\quad \qquad
\end{equation}
in the sense of isomorpy of normed spaces and
\begin{equation}\centering\label{3.4}
\begin{aligned}
\mathcal{P}\;\cong\;(\,\mathcal{L}^{2}\,, 
\, \|...\|^{\sim}\,) &\quad,\\ 
\text{e.g. with norm} \quad
&\|\,(\,a\,,\,b\,)\,\|^{\sim}\;:\,=\;
\sqrt{\| a \|^{2}\,+\,\| b \|^{2}}
\quad ;\quad a,b \in\mathcal{L}\;
\end{aligned}
\end{equation}
(or any equivalent norm).

As is well known, DSC and TLM algorithms follow a two-step iteration cycle
in working with alternate application of a
\emph{connection} and \emph{reflection} map
\vspace{-.2cm}
\begin{equation}\centering\notag %\label{3.5}
\mathcal{C}\,:\,I\times\mathcal{P}_{out}^{I}\,
\to\,\mathcal{P}_{in}
\qquad\text{and}\qquad 
\mathcal{R}\,:\,J\times\mathcal{P}_{in}^{J}\,
\to\,\mathcal{P}_{out} \quad,
\end{equation}
which respectively update the propagating fields at even and odd
integer multiples of half the time step, i.e. on ${I\,:\,=\,}$ ${
\{\,k\tau\mid k\in\mathbb{Z}\,\}}$ and
${J\,:\,=\,}$ ${\{\,(2k+1)\tau/2\mid k\in\mathbb{Z}\,\}}$.
To these maps the following functions $\,F_{_{C}}\,,
F_{_{R}}\,$ are associated in a one-to-one correspondence
\begin{equation}\centering\label{3.6}
\begin{aligned}
F_{_{C}}\,:\,\mathcal{L}^{I}\,&\to\,\mathcal{L}^{I} \\[-6pt]
f\,&\mapsto\,F_{_{C}}\,f\quad\text{with}\quad
F_{_{C}}\,f\,(\,t\,)\;:\,=\;nb\circ\mathcal{C}\,(\,t\,,\,f\,)
\quad
\end{aligned}
\end{equation}
and 
\begin{equation}\centering\label{3.7}
\begin{aligned}
F_{_{R}}\,:\,\mathcal{L}^{J}\,&\to\,\mathcal{L}^{J} \\[-6pt]
g\,&\mapsto\,F_{_{R}}\,g\quad\text{with}\quad
F_{_{R}}\,g\,(\,t\,)\;:\,=\;nb\circ\mathcal{R}\,(\,t\,,\,g\,)
\quad.
\end{aligned}
\end{equation}
\vspace{-.2cm}
\begin{lsproposition}\label{P.3.1}
\begin{itemize}
\item[]
\item[(i)]
$F_{_{C}}\,$ and $\,F_{_{R}}\,$ are causal on
$\mathcal{L}^{I}$ and $\mathcal{L}^{J}$, respectively.
\item[(ii)]
For $\,r,s\in J\,$ and $\,T_{s}:f(t)\,\mapsto\,f(t+s)\,$ the
\emph{shift operator} on $\,\mathcal{L}^{I\;\cup\,J}\,$,
\begin{equation}\centering\notag %\label{3.8}
(\quad r,s \leqslant 0\quad\text{and}\quad
\begin{cases}
F_{_{C}}\\[-2pt]
F_{_{R}}
\end{cases}
\!\!\!\!\text{is $\,\alpha$-passive on}\quad
\begin{cases}
\mathcal{L}^{I}\\[-2pt]
\mathcal{L}^{J}
\end{cases}
\!)\;\Longrightarrow\;
\end{equation}
\vspace{-5pt}
\begin{equation}\centering\notag
\;\,\Longrightarrow\;(\quad
\begin{cases}
\,T_{r}\circ F_{_{C}}\circ T_{s}\\[-2pt]
\,T_{r}\circ F_{_{R}}\circ T_{s}
\end{cases}
\!\!\!\!\text{is $\,\alpha$-passive on}\quad
\begin{cases}
\mathcal{L}^{J}\\[-2pt]
\mathcal{L}^{I}
\end{cases}
\!)\quad .
\end{equation}
\end{itemize}
\end{lsproposition}

\vspace{-4pt}\hspace{-.55cm}
Indeed, the first statement is only a trivial consequence of the definition
of the DSC propagators $\,\mathcal{R}\,$ and $\,\mathcal{C}\,$ as
functions on back-in-time running sequences of incident and outgoing states,
cf. ~\cite{He1}, section ~3. Then (ii) follows from Definition ~\ref{D.2.3}.
\hfill $\Box$

\begin{lsdefinition} %\label{D.3.1}
\begin{equation}\centering\notag
\text{The}\;
\begin{cases}
\text{connection map $\mathcal{C}$}\\[-2pt]
\text{reflection map $\mathcal{R}$}
\end{cases}
\text{is called $\,\alpha$-passive,\;iff}\quad
\begin{cases}
F_{_{C}}\\[-2pt]
F_{_{R}}
\end{cases}
\text{is $\,\alpha$-passive}\quad
\end{equation}
in the sense of Definition ~\ref{D.2.3}.
\end{lsdefinition}

\hspace{-.55cm} We now claim the main statement:
\begin{lstheorem}{\label{T.3.1}}
With every time limited excitation, the DSC process generated by
$\,\alpha$-passive reflection and connection maps is uniformly bounded,
hence stable.
\end{lstheorem}

\begin{proof}
It is sufficient to show that every finitely excited DSC process
that is generated by $\,\alpha$-passive $\mathcal{R}\,$ and $\mathcal{C}\,$
can be written as a pair of processes, either of which satisfies
Theorem ~\ref{T.2.1}.

With $\,H\,:\,=\,I\cup J\,=\,{\{k \tau/2\,|\,k\in\mathbb{Z}\}}\,$
and the measure on $\,H\,$ inherited from $\,I\,$ and $\,J\,$ jointly,
the space of all DSC processes is
\begin{equation}\centering\notag %\label{3.9}
\begin{aligned}
\mathcal{E}\;&:\,=\;\{\;h=(\,f,\,g\,)\in (\mathcal{L}^{2})^{\,H}\quad |
\quad f((2k+1)\,\tau/2\,)\,=\,f(\,k\tau\,) \\
\quad &\text{and}\quad
g(\,k\tau\,)\,=\,g((2k-1)\,\tau/2\,)\;,\quad
\text{for every $\,k\in\mathbb{Z}$}\;\}\;,
\end{aligned}
\end{equation}
i.e. the functions $\,f:H\to\mathcal{L}\,$ and $\,g:H\to\mathcal{L}\,$
in ${\,h\,=\,(\,f\,,\,g\,)\,\in\mathcal{E}\,}$ 
'switch' at even and odd integer multiples of $\,\tau/2\,$, respectively.
So, there is a natural bijection
\vspace{-0.2cm}
\begin{equation}\centering\notag %\label{3.10}
\begin{aligned}
\mathcal{E}\;
&\to\;\mathcal{L}^{I}\times\mathcal{L}^{J}\\[-3pt]
h\,=\,(\,f\,,\,g\,)\;
&\mapsto\;(\,f\downharpoonright I\,,\,g\downharpoonright J\,)\quad,
\end{aligned}
\end{equation}
in virtue of which the first and second components in
$\,\mathcal{E}\,$ can be naturally identified with $\,\mathcal{L}^{I}\,$
and $\,\mathcal{L}^{J}\,$, respectively.

For any incident function
$\,e\,\in\,\mathcal{L}^{I}\times\{\,0\,\}\,\subset\mathcal{E}\,$
supported on a finite interval $\,[\,0\,,\,N\tau\,)\subset I\,$,
the DSC process excited with $\,e\,$ and generated by
$\,\mathcal{C}\,$ and $\,\mathcal{R}\,$ is the well-defined function
$\,h\in\,(\,\mathcal{L}^{2}\,)^{H}\,$ 
such that $\,h\,(\,t\,)\,=\,0\,$ for $\,t\,\leqslant\,0\,$ and recursively
for $\,0\,<\,t\,\in\,H\,$
\vspace{-.2cm}
\begin{equation}\centering
\label{3.11}
\begin{aligned}
h\,(\,t\,+\,\frac{\tau}{2}\,)\,=\,
\begin{cases}
(\,h_{1}\,(\,t\,)\,,\;T_{-\frac{\tau}{2}}\,F_{_{R}}\,
T_{-\frac{\tau}{2}}\,[\,e\,+\,h_{1}\,]\,(\,t\,)\,)
&\,\text{if $\,t\,\in\,I\,$}\\
(\,T_{-\frac{\tau}{2}}\,F_{_{C}}\,
T_{-\frac{\tau}{2}}\,[\,h_{2}\,]\,(\,t\,)\,,\;h_{2}(\,t\,)\,)
&\,\text{if $\,t\,\in\,J\,$}\,.
\end{cases}
\end{aligned}
\end{equation}
Actually $\,h\,$ is well and uniquely defined by relations \eqref{3.11} which
provide separate recurrence relations for the two processes
$\,h_{1}\,,\,h_{2}\,$;
for instance for $\,h_{1}\,$ and $\,t\,\in I\,$
\vspace{-4pt}
\begin{equation}\centering
\label{3.12}
\begin{aligned}
h_{1}\,(\,t\,+\,\tau\,)\;&=\quad h_{1}\,(\,t^{\sptilde}\!+\,\frac{\tau}{2}\,)
\quad\text{with}\quad t^{\sptilde}:\,=\,t\,+\,\frac{\tau}{2}\,\in\,J\\
&=\quad T_{\frac{-\tau}{2}}F_{_{C}}\,T_{\frac{-\tau}{2}}
[\,h_{2}\mid _{t^{\sptilde}}\,]\,(\,t^{\sptilde})
&&\text{by \eqref{3.11}}\\
&=\quad F_{_{C}}T_{\frac{-\tau}{2}}[\hspace{-58pt}
\underbrace{h_{2}\mid _{t^{\sptilde}}}_{\hspace{76pt}=\;T_{\frac{\tau}{2}} h_{2}
\mid _{t^{\sptilde}}\;\text{,\;cf.\;\eqref{3.11}}}
\hspace{-58pt}]\,(\,t^{\sptilde}\!-\frac{\tau}{2}\,)
&&\,\\
&=\quad F_{_{C}}T_{\frac{-\tau}{2}}\,F_{_{R}}\,T_{\frac{-\tau}{2}}
[\,h_{1}+\,e\,]\,(\,t\,)\quad .
\end{aligned}
\end{equation}
In virtue of the causality of
$F_{_{C}}T_{\frac{-\tau}{2}}\,F_{_{R}}\,T_{\frac{-\tau}{2}}$\,
function $\,h_{1}\,$ enters the last \linebreak expression
in \eqref{3.12} only up to argument $\,t\,$.
Moreover, since products of $\,\alpha$-passive operators are
obviously $\alpha$-passive,
${F_{_{C}}\,T_{-\tau/2}\,F_{_{R}}\,T_{-\tau/2}}$ is $\,\alpha$-passive
due to Proposition ~\ref{P.3.1} ~(ii). Hence Theorem~\ref{T.2.1}
applies to $\,h_{1}\,$, which thus is stable, just as then also is
$\,h_{2}\,=\,F_{_{R}}\,T_{-\tau/2}\,[\,e+h_{1}\,]\,$.
\end{proof}

%% file: ccl.tex
\section{Conclusions}\label{C:sec4}
It has been demonstrated that a class of $\,\alpha$-passive processes,
characterised by simple contraction properties with respect to a
\emph{limiting} functional $\,\alpha\,$, are unconditionally stable.
Finitely excited DSC processes, generated by \linebreak
$\,\alpha$-passive reflection and connection maps, are then necessarily
stable. \linebreak
Stability analysis reduces so in pursuant cases to finding \emph{any}
limiting functional $\,\alpha\,$ in relation to which the
reflection and connection maps are \linebreak
$\,\alpha$-passive. In reverse: If the connection and reflection maps
are \emph{per se} so designed as to conserve (\,or admit loss of\;) any
quantity 'measured' by a limiting functional, then the DSC process generated
by those maps is essentially stable.

%% file: ref.tex
%\newpage

%% file: stb.bbl
\begin{thebibliography}{DGMS}

\bibitem[JoBe]{JoBe}
Johns, P.B. and Beurle, R.L. Numerical solution of 2-dimensional scattering
problems using transmission line matrix, Proc. IEEE, vol. 118, pp. 1203-1208,
1971

\bibitem[Jo1]{Jo1}
Akhtarzad, S., Johns, P.B., Solution of 6-components electromagnetic fields
in three dimensions and time by the T.L.M. method, Electron. Lett., vol. 10,
pp. 535-537, Dec. 12, 1974

\bibitem[Jo2]{Jo2}
Johns, P.B., A simple, explicit and unconditionally stable numerical routine
for the solution of the diffusion equations, Int. J. Numer. Methods Eng.,
vol. 11, pp. 1307-1328, 1977 

\bibitem[Ch]{Ch}
Christopoulos, C., {\it The Transmission-Line Modeling Method TLM},
IEEE Press, New York 1995

\bibitem[dC]{dC}
De Cogan, D., {\it Transmission Line Matrix (TLM) Techniques for Diffusion
Applications}, Gordon and Breach, 1998

\bibitem[Tlm1]{Tlm1}
Proceedings of the 1st Int. Workshop on Transmission Line Matrix (TLM)
Modelling, Victoria, 1995

\bibitem[Tlm2]{Tlm2}
Proceedings of the 2nd Int. Workshop on Transmission Line Matrix (TLM)
Modelling, TU M\"unchen, 1997

\bibitem[Tlm3]{Tlm3}
Proceedings of the 3nd Int. Workshop on Transmission Line Matrix (TLM)
Modelling, Nice, 1999

\bibitem[Re]{Re}
Rebel, J.N., On the Foundations of the Transmission Line Matrix Method,
Thesis, TU M\"unchen, 2000

\bibitem[Hoe]{Hoe}
Hoefer, W. J. R., The Transmission Line Matrix (TLM) Method, in Numerical
Techniques for Microwave and Millimeter-Wave Passive Structures,
T. Itoh ed. John Wiley \& Sons, 1989

\bibitem[Yee1]{Yee1}
Yee, K.S., Numerical Solution of Initial Boundary Value Problems Invoking
Maxwell's Equations in Isotropic Media, IEEE Trans. Antennas and Propagation,
vol. 14 (3), pp 302-307, May 1966

\bibitem[He1]{He1}
\makebox[80pt][l]{Hein, S., Dual scattering channel schemes extending the
Johns Algorithm,} \linebreak http://arxiv.org/abs/math/NA0309261, Sep 2003
(preprint)

\bibitem[He2]{He2}
Hein, S., Finite-difference time-domain approximation of Maxwell's equations
with nonorthogonal condensed TLM mesh, Int. J. Num. Modelling, 
vol. 7, pp. 179-188, 1994 

\bibitem[He3]{He3}
Hein, S., Synthesis of TLM Algorithms in the Popagator Integral Framework,
Proceedings of the 2nd. Int. Workshop on Transmission  Line Matrix Modeling
(TLM) - Theory and Applications, pp. 1-11, Munich, October 1997 (invited paper)

\bibitem[He4]{He4}
Hein, S., TLM numerical solution of Bloch's equations for magnetized gyrotropic
media, Appl. Math. Modelling, vol. 21, pp. 221-229, 1997

\bibitem[He5]{He5}
Hein, S., A TLM node for superconducting boundary illustrating the
propagator approach, Spinner Report E017, M\"unchen 1992

\bibitem[He6]{He6}
Hein, S., Consistent finite difference modelling of Maxwell's equations
with lossy symmetrical condensed TLM node,  Int. J. Num. Modelling,
vol. 6, pp. 207-220, 1993

\end{thebibliography}
